# Stability of Matrix Factorization for Collaborative Filtering

**Yu-Xiang Wang**                                              YUXIANGWANG@NUS.EDU.SG

Department of Mechanical Engineering, National University of Singapore, Singapore 117576

**Huan Xu**                                                      MPEXUH@NUS.EDU.SG

Department of Mechanical Engineering, National University of Singapore, Singapore 117576

## Abstract

We study the stability *vis a vis* adversarial noise of matrix factorization algorithm for matrix completion. In particular, our results include: (I) we bound the gap between the solution matrix of the factorization method and the ground truth in terms of root mean square error; (II) we treat the matrix factorization as a subspace fitting problem and analyze the difference between the solution subspace and the ground truth; (III) we analyze the prediction error of individual users based on the subspace stability. We apply these results to the problem of collaborative filtering under manipulator attack, which leads to useful insights and guidelines for collaborative filtering system design.

## 1. Introduction

Collaborative prediction of user preferences has attracted fast growing attention in the machine learning community, best demonstrated by the million-dollar Netflix Challenge. Among various models proposed, matrix factorization is arguably the most widely applied method, due to its high accuracy, scalability (Su & Khoshgoftaar, 2009) and flexibility to incorporating domain knowledge (Koren et al., 2009). Hence, not surprisingly, matrix factorization is the centerpiece of most state-of-the-art collaborative filtering systems, including the winner of Netflix Prize (Bell & Koren, 2007). Indeed, matrix factorization has been widely applied to tasks other than collaborative filtering, including structure from motion, localization in wireless sensor network, DNA microarray estimation and beyond. Matrix factorization is also considered as a fundamental building block of many popular algorithms in regression, factor analysis, dimension reduction, and clustering (Singh & Gordon, 2008).

Despite the popularity of factorization methods, not much has been done on the theoretical front. In this paper, we fill the blank by analyzing the stability *vis a vis* adversarial noise of the matrix factorization methods, in hope of providing useful insights and guidelines for practitioners to design and diagnose their algorithm efficiently.

Our main contributions are three-fold: In Section 3 we bound the gap between the solution matrix of the factorization method and the ground truth in terms of root mean square error. In Section 4, we treat the matrix factorization as a subspace fitting problem and analyze the difference between the solution subspace and the ground truth. This facilitates an analysis of the prediction error of individual users, which we present in Section 5. To validate these results, we apply them to the problem of collaborative filtering under manipulator attack in Section 6. Interestingly, we find that matrix factorization are robust to the so-called "targeted attack", but not so to the so-called "mass attack" unless the number of manipulators are small. These results agree with the simulation observations.

We briefly discuss relevant literatures. Azar et al. (2001) analyzed asymptotic performance of matrix factorization methods, yet under stringent assumptions on the fraction of observation and on the singular values. Drineas et al. (2002) relaxed these assumptions but it requires a few fully rated users – a situation that rarely happens in practice. Srebro (2004) considered the problem of the generalization error of learning a low-rank matrix. Their technique is similar to the proof of our first result, yet applied to a different context. Specifically, they are mainly interested in *binary* prediction (i.e., "like/dislike") rather than recovering the *real-valued* ground-truth matrix (and its column subspace). In addition, they did not investigate the stability of the algorithm under noise and manipulators.

---

Recently, some alternative algorithms, notably StableMC (Candes & Plan, 2010) based on nuclear norm optimization, and OptSpace (Keshavan et al., 2010b) based on gradient descent over the Grassmannian, have been shown to be stable *vis a vis* noise (Candes & Plan, 2010; Keshavan et al., 2010a). However, these two methods are less effective in practice. As documented in Mitra et al. (2010); Wen (2010) and many others, matrix factorization methods typically outperform these two methods. Indeed, our theoretical results reassure these empirical observations, see Section 3 for a detailed comparison of the stability results of different algorithms.

## 2. Formulation

### 2.1. Matrix Factorization with Missing Data

Let the user ratings of items (such as movies) form a matrix $Y$, where each column corresponds to a user and each row corresponds to an item. Thus, the $ij^{th}$ entry is the rating of item-$i$ from user-$j$. The valid range of the rating is $[-k, +k]$. $Y$ is assumed to be a rank-$r$ matrix[1], so there exists a factorization of this rating matrix $Y = UV^T$, where $Y \in \mathbb{R}^{m \times n}$, $U \in \mathbb{R}^{m \times r}$, $V \in \mathbb{R}^{n \times r}$. Without loss of generality, we assume $m \le n$ throughout the paper.

Collaborative filtering is about to recover the rating matrix from a fraction of entries possibly corrupted by noise or error. That is, we observe $\widehat{Y}_{ij}$ for $(ij) \in \Omega$ the sampling set (assumed to be uniformly random), and $\widehat{Y} = Y + E$ being a corrupted copy of $Y$, and we want to recover $Y$. This naturally leads to the optimization program below:

$$\min_{U,V} \quad \frac{1}{2} \left\| P_\Omega(UV^T - \widehat{Y}) \right\|_F^2$$
$$\text{subject to} \quad \left| [UV^T]_{i,j} \right| \le k, \tag{1}$$

where $P_\Omega$ is the sampling operator defined to be:

$$[P_\Omega(Y)]_{i,j} = \begin{cases} Y_{i,j} & \text{if } (i,j) \in \Omega; \\ 0 & \text{otherwise.} \end{cases} \tag{2}$$

We denote the optimal solution $Y^* = U^*V^{*T}$ and the error $\Delta = Y^* - Y$.

### 2.2. Matrix Factorization as Subspace Fitting

As pointed out in Chen (2008), an alternative interpretation of collaborative filtering is fitting the optimal $r$-dimensional subspace $\mathcal{N}$ to the sampled data. That

is, one can reformulate (1) into an equivalent form[2]:

$$\min_N f(N) = \sum_i \|(I - \mathbb{P}_i)y_i\|^2 = \sum_i y_i^T(I - \mathbb{P}_i)y_i, \tag{3}$$

where $y_i$ is the observed entries in the $i^{th}$ column of $Y$, $N$ is an $m \times r$ matrix representing an orthonormal basis[3] of $\mathcal{N}$, $N_i$ is the restriction of $N$ to the observed entries in column $i$, and $\mathbb{P}_i = N_i(N_i^T N_i)^{-1}N_i^T$ is the projection onto span($N_i$).

After solving (3), we can estimate the full matrix in a column by column manner via (4). Here $y_i^*$ denotes the full $i^{th}$ column of recovered rank-$r$ matrix $Y^*$.

$$y_i^* = N(N_i^T N_i)^{-1}N_i^T y_i = N\text{pinv}(N_i)y_i. \tag{4}$$

Due to error term $E$, the ground truth subspace $\mathcal{N}^{gnd}$ can not be obtained. Instead, denote the optimal subspace of (1) (equivalently (3)) by $\mathcal{N}^*$, and we bound the gap between these two subspaces using Canonical angle. The canonical angle matrix $\Theta$ is an $r \times r$ diagonal matrix, with the $i^{th}$ diagonal entry $\theta_i = \arccos \sigma_i((N^{gnd})^T N^*)$.

The error of subspace recovery is measured by $\rho = \|\sin \Theta\|_2$, justified by the following properties adapted from Chapter 2 of Stewart & Sun (1990):

$$\|\mathbb{P}^{gnd} - \mathbb{P}^{\mathcal{N}^*}\|_F = \sqrt{2}\|\sin \Theta\|_F,$$
$$\|\mathbb{P}^{gnd} - \mathbb{P}^{\mathcal{N}^*}\|_2 = \|\sin \Theta\|_2 = \sin \theta_1. \tag{5}$$

### 2.3. Algorithms

We focus on the stability of the *global optimal solution* of Problem (1). As Problem (1) is not convex, finding the global optimum is non-trivial in general. While this is certainly an important question, it is beyond the scope of this paper. Instead, we briefly review some results on this aspect.

The simplest algorithm for (1) is perhaps the alternating least square method (ALS) which alternately minimizes the objective function over $U$ and $V$ until convergence. More sophisticatedly, second-order algorithms such as Wiberg, Damped Newton and Levenberg Marquadt are proposed with better convergence rate, as surveyed in Okatani & Deguchi (2007). Specific variations for CF are investigated in Takács et al. (2008) and Koren et al. (2009).

From an empirical perspective, Mitra et al. (2010) reported that the global optimum is often obtained in simulation and Chen (2008) demonstrated satisfactory percentage of hits to global minimum from randomly initialized trials on a real data set.

---

[1] In practice, this means the user's preference of movies are influenced by no more than $r$ latent factors.

[2] Strictly speaking, this is only equivalent to (1) without the box constraint. See the discussion in Supplementary Material for our justifications.

[3] It is easy to see $N = \text{ortho}(U)$ for $U$ in (1)

## 3. Stability

We show in this section that when *sufficiently many* entries are sampled, the global optimal solution of factorization methods is stable *via a vis* noise – i.e., it recovers a matrix "close to" the ground-truth. This is measured by the root mean square error (RMSE):

$$\text{RMSE} = \frac{1}{\sqrt{mn}}\|Y^* - Y\| \qquad (6)$$

**Theorem 1.** *There exists an absolute constant $C$, such that with probability at least $1 - 2\exp(-n)$,*

$$\text{RMSE} \leq \frac{1}{\sqrt{|\Omega|}}\|P_\Omega(E)\|_F + \frac{\|E\|_F}{\sqrt{mn}} + Ck\left(\frac{nr\log(n)}{|\Omega|}\right)^{\frac{1}{4}}.$$

Notice that when $|\Omega| \gg nr\log(n)$ the last term diminishes, and the RMSE is essentially bounded by the "average" magnitude of entries of $E$, i.e., the factorization methods are stable.

### Comparison with related work

We recall similar RMSE bounds for StableMC of Candes & Plan (2010) and OptSpace of Keshavan et al. (2010a):

**StableMC:** RMSE
$$\leq \sqrt{\frac{32\min(m,n)}{|\Omega|}}\|P_\Omega(E)\|_F + \frac{1}{\sqrt{mn}}\|P_\Omega(E)\|_F. \qquad (7)$$

**OptSpace:** $\text{RMSE} \leq C\kappa^2\frac{n\sqrt{r}}{|\Omega|}\|P_\Omega(E)\|_2. \qquad (8)$

Albeit the fact that these bounds are for different algorithms and under different assumptions (see Table 1 for details), it is still interesting to compare the results with Theorem 1. We observe that Theorem 1 is tighter than (7) by a scale of $\sqrt{\min(m,n)}$, and tighter than (8) by a scale of $\sqrt{n/\log(n)}$ in case of adversarial noise. However, the latter result is stronger when the noise is stochastic, due to the spectral norm used.

### Compare with an Oracle

We next compare the bound with an oracle, introduced in Candes & Plan (2010), that is assumed to know the ground-truth column space $\mathcal{N}$ *a priori* and recover the matrix by projecting the observation to $\mathcal{N}$ in the least square sense column by column via (4). It is shown that RMSE of this oracle satsifies,

$$\text{RMSE} \approx \sqrt{1/|\Omega|}\|P_\Omega(E)\|_F. \qquad (9)$$

Notice that Theorem 1 matches this oracle bound, and hence it is tight up to a constant factor.

### 3.1. Proof of Stability Theorem

We briefly explain the proof idea first. By definition, the algorithm finds the optimal rank-$r$ matrices, measured in terms of the root mean square (RMS) on the *sampled* entries. To show this implies a small RMS on the *entire* matrix, we need to bound their gap

$$\tau(\Omega) \triangleq \Big|\frac{1}{\sqrt{|\Omega|}}\|P_\Omega(\widehat{Y} - Y^*)\|_F - \frac{1}{\sqrt{mn}}\|\widehat{Y} - Y^*\|_F\Big|.$$

To bound $\tau(\Omega)$, we require the following theorem.

**Theorem 2.** *Let $\hat{\mathcal{L}}(X) = \frac{1}{\sqrt{|\Omega|}}\|P_\Omega(X - \widehat{Y})\|_F$ and $\mathcal{L}(X) = \frac{1}{\sqrt{mn}}\|X - \widehat{Y}\|_F$ be the empirical and actual loss function respectively. Furthermore, assume entry-wise constraint $\max_{i,j}|X_{i,j}| \leq k$. Then for all rank-$r$ matrices $X$, with probability greater than $1 - 2\exp(-n)$, there exists a fixed constant $C$ such that*

$$\sup_{X\in S_r}|\hat{\mathcal{L}}(X) - \mathcal{L}(X)| \leq Ck\left(\frac{nr\log(n)}{|\Omega|}\right)^{\frac{1}{4}}.$$

Indeed, Theorem 2 easily implies Theorem 1.

*Proof of Theorem 1.* The proof makes use of the fact that $Y^*$ is the global optimal of (1).

$$\begin{aligned}
\text{RMSE} &= \frac{1}{\sqrt{mn}}\|Y^* - Y\|_F = \frac{1}{\sqrt{mn}}\|Y^* - \widehat{Y} + E\|_F \\
&\leq \frac{1}{\sqrt{mn}}\|Y^* - \widehat{Y}\|_F + \frac{1}{\sqrt{mn}}\|E\|_F \\
&\overset{(a)}{\leq} \frac{1}{\sqrt{|\Omega|}}\|P_\Omega(Y^* - \widehat{Y})\|_F + \tau(\Omega) + \frac{1}{\sqrt{mn}}\|E\|_F \\
&\overset{(b)}{\leq} \frac{1}{\sqrt{|\Omega|}}\|P_\Omega(Y - \widehat{Y})\|_F + \tau(\Omega) + \frac{1}{\sqrt{mn}}\|E\|_F \\
&= \frac{1}{\sqrt{|\Omega|}}\|P_\Omega(E)\|_F + \tau(\Omega) + \frac{1}{\sqrt{mn}}\|E\|_F.
\end{aligned}$$

Here, (a) holds from definition of $\tau(\Omega)$, and (b) holds because $Y^*$ is optimal solution of (1). Since $Y^* \in S_r$, applying Theorem 2 completes the proof. $\qquad\square$

The proof of Theorem 2 is deferred to Appendix A due to space constraints. The main idea, briefly speaking, is to bound, for a fixed $X \in S_r$,

$$\begin{aligned}
&\big|(\hat{\mathcal{L}}(X))^2 - (\mathcal{L}(X))^2\big| \\
&= \big|\frac{1}{|\Omega|}\|P_\Omega(X - \widehat{Y})\|_F^2 - \frac{1}{mn}\|X - \widehat{Y}\|_F^2\big|,
\end{aligned}$$

using Hoeffding's inequality for sampling without replacement; then bound $\big|\hat{\mathcal{L}}(X) - \mathcal{L}(X)\big|$ using

$$\big|\hat{\mathcal{L}}(X) - \mathcal{L}(X)\big| \leq \sqrt{\big|(\hat{\mathcal{L}}(X))^2 - (\mathcal{L}(X))^2\big|};$$

and finally, bound $\sup_{X\in S_r}|\hat{\mathcal{L}}(X) - \mathcal{L}(X)|$ using an $\epsilon$−net argument.

| | Rank constraint | $Y_{i,j}$ constraint | $\sigma$ constraint | incoherence | global optimal |
|---|---|---|---|---|---|
| Theorem 1 | fixed rank | box constraint | no | no | assumed |
| OptSpace | fixed rank | regularization | condition number | weak | not necessary |
| NoisyMC | Relaxed to trace | implicit | no | strong | yes |

*Table 1.* Comparison of assumptions between stability results in our Theorem 1, OptSpace and NoisyMC

## 4. Subspace Stability

In this section we investigate the stability of recovered *subspace* using matrix factorization methods. Recall that matrix factorization methods assume that, in the idealized noiseless case, the preference of each user belongs to a low-rank subspace. Therefore, if this subspace can be readily recovered, then we can predict preferences of a new user without re-run the matrix factorization algorithms. We analyze the latter, prediction error on individual users, in Section 5.

To illustrate the difference between the stability of the recovered matrix and that of the recovered subspace, consider a concrete example in movie recommendation, where there are both honest users and malicious manipulators in the system. Suppose we obtain an output subspace $N^*$ by (3) and the missing ratings are filled in by (4). If $N^*$ is very "close" to ground truth subspace $N$, then all the predicted ratings for honest users will be good. On the other hand, the prediction error of the preference of the manipulators – who do not follow the low-rank assumption – can be large, which leads to a large error of the recovered matrix. Notice that we are only interested in predicting the preference of the honest users. Hence the subspace stability provides a more meaningful metric here.

### 4.1. Subspace Stability Theorem

Let $\mathcal{N}, \mathcal{M}$ and $\mathcal{N}^*, \mathcal{M}^*$ be the $r$-dimensional column space-row space pair of matrix $Y$ and $Y^*$ respectively. We'll denote the corresponding $m \times r$ and $n \times r$ orthonormal basis matrix of the vector spaces using $N, M, N^*, M^*$. Furthermore, Let $\Theta$ and $\Phi$ denote the canonical angles $\angle(\mathcal{N}^*, \mathcal{N})$ and $\angle(\mathcal{M}^*, \mathcal{M})$ respectively.

**Theorem 3.** *When $Y$ is perturbed by additive error $E$ and observed only on $\Omega$, then there exists a $\Delta$ satisfying* $\|\Delta\| \leq \sqrt{\frac{mn}{|\Omega|}} \|P_\Omega(E)\|_F + \|E\|_F + \sqrt{mn} \, |\tau(\Omega)|$, *such that:*

$$\|\sin\Theta\| \leq \frac{\sqrt{2}}{\delta} \|(\mathbb{P}^{\mathcal{N}^\perp}\Delta)\|; \quad \|\sin\Phi\| \leq \frac{\sqrt{2}}{\delta} \|(\mathbb{P}^{\mathcal{M}^\perp}\Delta^T)\|,$$

*where $\|\cdot\|$ is either the Frobenious norm or the spectral norm, and $\delta = \sigma_r^*$, i.e., the $r^{th}$ largest singular value of the recovered matrix $Y^*$.*

*Furthermore, we can bound $\delta$ by:*

$$
\begin{cases}
\sigma_r - \|\Delta\|_2 & \leq \delta \leq & \sigma_r + \|\Delta\|_2 \\
\sigma_r^{\tilde{Y}_\mathcal{N}} - \|\mathbb{P}^{\mathcal{N}^\perp}\Delta\|_2 & \leq \delta \leq & \sigma_r^{\tilde{Y}_\mathcal{N}} + \|\mathbb{P}^{\mathcal{N}^\perp}\Delta\|_2 \\
\sigma_r^{\tilde{Y}_\mathcal{M}} - \|\mathbb{P}^{\mathcal{M}^\perp}\Delta^T\|_2 & \leq \delta \leq & \sigma_r^{\tilde{Y}_\mathcal{M}} + \|\mathbb{P}^{\mathcal{M}^\perp}\Delta^T\|_2
\end{cases}
$$

*where $\tilde{Y}_\mathcal{N} = Y + \mathbb{P}^\mathcal{N}\Delta$ and $\tilde{Y}_\mathcal{M} = Y + (\mathbb{P}^\mathcal{M}\Delta^T)^T$.*

Notice that in practice, as $Y^*$ is the output of the algorithm, its $r^{th}$ singular value $\delta$ is readily obtainable. Intuitively, Theorem 3 shows that the subspace sensitivity *vis a vis* noise depends on the singular value distribution of original matrix $Y$. A well-conditioned rank-$r$ matrix $Y$ can tolerate larger noise, as its $r^{th}$ singular value is of the similar scale to $\|Y\|_2$, its largest singular value.

### 4.2. Proof of Subspace Stability

*Proof of Theorem 3.* In the proof, we use $\|\cdot\|$ when a result holds for both Frobenious norm and for spectral norm. We prove the two parts separately.

*Part 1: Canonical Angles.*

Let $\Delta = Y^* - Y$. By Theorem 1, we have $\|\Delta\| \leq \sqrt{\frac{mn}{|\Omega|}} \|P_\Omega(E)\|_F + \|E\|_F + \sqrt{mn} \, |\tau(\Omega)|$. The rest of the proof relates $\Delta$ with the deviation of spaces spanned by the top $r$ singular vectors of $Y$ and $Y^*$ respectively. Our main tools are Weyl's Theorem and Wedin's Theorem (Lemma F.1 and F.2 in Appendix F).

We express singular value decomposition of $Y$ and $Y^*$ in block matrix form as in (F.1) and (F.2) of Appendix F, and set the dimension of $\Sigma_1$ and $\hat{\Sigma}_1$ to be $r \times r$. Recall, rank$(Y) = r$, so $\Sigma_1 = \text{diag}(\sigma_1, ..., \sigma_r)$, $\Sigma_2 = 0$, $\hat{\Sigma}_1 = \text{diag}(\sigma_1', ..., \sigma_r')$. By setting $\hat{\Sigma}_2$ to 0 we obtained $Y'$, the nearest rank-$r$ matrix to $Y^*$. Observe that $N^* = \hat{L}_1$, $M^* = (\hat{R}_1)^T$.

To apply Wedin's Theorem (Lemma F.2), we have the residual $Z$ and $S$ as follows:

$$
\begin{aligned}
Z &= YM^* - N^*\hat{\Sigma}_1, \\
S &= Y^T N^* - M^*\hat{\Sigma}_1,
\end{aligned}
$$

which leads to

$$
\begin{aligned}
\|Z\| &= \|(\hat{Y} - \Delta)M^* - N^*\hat{\Sigma}_1\| = \|\Delta M^*\|, \\
\|S\| &= \|(\hat{Y} - \Delta)^T N^* - M^*\hat{\Sigma}_1\| = \|\Delta^T N^*\|.
\end{aligned}
$$

Substitute this into the Wedin's inequality, we have

$$\sqrt{\|\sin\Phi\|^2 + \|\sin\Theta\|^2} \leq \frac{\sqrt{\|\Delta^T N'\|^2 + \|\Delta M'\|^2}}{\delta}, \tag{10}$$

where $\delta$ satisfies (F.3) and (F.4). Specifically, $\delta = \sigma_r^*$. Observe that Equation (10) implies

$$\|\sin\Theta\| \leq \frac{\sqrt{2}}{\delta}\|\Delta\|; \quad \|\sin\Phi\| \leq \frac{\sqrt{2}}{\delta}\|\Delta\|.$$

To reach the equations presented in the theorem, we can tighten the above bound by decomposing $\Delta$ into two orthogonal components.

$$Y^* = Y + \Delta = Y + \mathbb{P}^{\mathcal{N}}\Delta + \mathbb{P}^{\mathcal{N}^\perp}\Delta := \tilde{Y}^{\mathcal{N}} + \mathbb{P}^{\mathcal{N}^\perp}\Delta. \tag{11}$$

It is easy to see that column space of $Y$ and $\tilde{Y}_\mathcal{N}$ are identical. So the canonical angle $\Theta$ between $Y^*$ and $Y$ are the same as that between $Y^*$ and $\tilde{Y}_\mathcal{N}$. Therefore, we can replace $\Delta$ by $\mathbb{P}^{\mathcal{N}^\perp}\Delta$ to obtain the equation presented in the theorem. The corresponding result for row subspace follows similarly, by decomposing $\Delta^T$ to its projection on $\mathcal{M}$ and $\mathcal{M}^\perp$.

*Part 2: Bounding $\delta$.*

We now bound $\delta$, or equivalently $\sigma_r^*$. By Weyl's theorem (Lemma F.1), we have

$$|\delta - \sigma_r| < \|\Delta\|_2.$$

Moreover, Applying Weyl's theorem on Equation (11), we have

$$|\delta - \sigma_r^{\tilde{Y}_\mathcal{N}}| \leq \|\mathbb{P}_{\mathcal{N}^\perp}\Delta\|_2.$$

Similarly, we have

$$|\delta - \sigma_r^{\tilde{Y}_\mathcal{M}}| \leq \|\mathbb{P}_{\mathcal{M}^\perp}\Delta^T\|_2.$$

This establishes the theorem. □

# 5. Prediction Error of individual user

In this section, we analyze how confident we can predict the ratings of a new user $y \in \mathcal{N}^{gnd}$, based on the subspace recovered via matrix factorization methods. In particular, we bound the prediction $\|\tilde{y}^* - y\|$, where $\tilde{y}^*$ is the estimation from partial rating using (4), and $y$ is the ground truth.

Without loss of generality, if the sampling rate is $p$, we assume observations occur in first $pm$ entries, such that $y = \begin{pmatrix} y_1 \\ y_2 \end{pmatrix}$ with $y_1$ observed and $y_2$ unknown.

## 5.1. Prediction of $y$ With Missing data

**Theorem 4.** *With all the notations and definitions above, and let $N_1$ denote the restriction of $N$ on the observed entries of $y$. Then the prediction for $y \in \mathcal{N}^{gnd}$ has bounded performance:*

$$\|\tilde{y}^* - y\| \leq \left(1 + \frac{1}{\sigma_{min}}\right)\rho\|y\|,$$

*where $\rho = \|\sin\Theta\|$ (see Theorem 3), $\sigma_{min}$ is the smallest non-zero singular value of $N_1$ ($r^{th}$ when $N_1$ is non-degenerate).*

*Proof.* By (4), and recall that only the first $pm$ entries are observed, we have

$$\tilde{y}^* = N \cdot \mathrm{pinv}(N_1)y_1 := \begin{pmatrix} y_1 - \tilde{e}_1 \\ y_2 - \tilde{e}_2 \end{pmatrix} := y + \tilde{e}.$$

Let $y^*$ be the vector obtained by projecting $y$ onto subspace $N$, and denote $y^* = \begin{pmatrix} y_1^* \\ y_2^* \end{pmatrix} = \begin{pmatrix} y_1 - e_1 \\ y_2 - e_2 \end{pmatrix} = y - e$, we have:

$$\begin{aligned} \tilde{y}^* &= N \cdot \mathrm{pinv}(N_1)(y_1^* + e_1) \\ &= N \cdot \mathrm{pinv}(N_1)y_1^* + N \cdot \mathrm{pinv}(N_1)e_1 \\ &= y^* + N \cdot \mathrm{pinv}(N_1)e_1. \end{aligned}$$

Then

$$\begin{aligned} \|\tilde{y}^* - y\| &= \|y^* - y + N \cdot \mathrm{pinv}(N_1)e_1\| \\ &\leq \|y^* - y\| + \frac{1}{\sigma_{min}}\|e_1\| \\ &\leq \rho\|y\| + \frac{1}{\sigma_{min}}\|e_1\|. \end{aligned}$$

Finally, we bound $e_1$ as follows

$$\|e_1\| \leq \|e\| = \|y - y^*\| \leq \|(\mathbb{P}^{gnd} - \mathbb{P}^{\mathcal{N}})y\| \leq \rho\|y\|,$$

which completes the proof. □

Suppose $y \notin \mathcal{N}^{gnd}$ and $y = \mathbb{P}^{gnd}y + (I - \mathbb{P}^{gnd})y := y^{gnd} + y^{gnd^\perp}$, then we have

$$\|e_1\| \leq \|(\mathbb{P}^{gnd} - \mathbb{P}^{\mathcal{N}})y\| + \|y^{gnd^\perp}\| \leq \rho\|y\| + \|y^{gnd^\perp}\|,$$

which leads to

$$\|\tilde{y}^* - y^{gnd}\| \leq \left(1 + \frac{1}{\sigma_{min}}\right)\rho\|y\| + \frac{\|y^{gnd^\perp}\|}{\sigma_{min}}.$$

## 5.2. Bound on $\sigma_{min}$

To complete the above analysis, we now bound $\sigma_{min}$. Notice that in general $\sigma_{min}$ can be arbitrarily close to zero, if $N$ is "spiky". Hence we impose the strong incoherence property introduced in Candes & Tao (2010) (see Appendix C for the definition) to avoid such situation. Due to space constraint, we defer the proof of the following to the Appendix C.

**Proposition 1.** *If matrix $Y$ satisfies strong incoherence property with parameter $\mu$, then:*

$$\sigma_{min}(N_1) \geq 1 - \left( \frac{r}{m} + (1-p)\mu\sqrt{r} \right)^{\frac{1}{2}}.$$

FOR GAUSSIAN RANDOM MATRIX

Stronger results on $\sigma_{min}$ is possible for randomly generated matrices. As an example, we consider the case that $Y = UV$ where $U$, $V$ are two Gaussian random matrices of size $m \times r$ and $r \times n$, and show that $\sigma_{min}(N_1) \approx \sqrt{p}$.

**Proposition 2.** *Let $G \in \mathbb{R}^{m \times r}$ have i.i.d. zero-mean Gaussian random entries. Let $N$ be its orthonormal basis[4]. Then there exists an absolute constant $C$ such that with probability of at least $1 - Cn^{-10}$,*

$$\sigma_{min}(N_1) \geq \sqrt{\frac{k}{m}} - 2\sqrt{\frac{r}{m}} - C\sqrt{\frac{\log m}{m}}.$$

Due to space limit, the proof of Proposition 2 is deferred to the Supplementary material. The main idea is to apply established results about the singular values of Gaussian random matrix $G$ (e.g., Rudelson & Vershynin, 2009; Silverstein, 1985; Davidson & Szarek, 2001), then show that the orthogonal basis $N$ of $G$ is very close to $G$ itself.

We remark that the bound on singular values we used has been generalized to random matrices following subgaussian (Rudelson & Vershynin, 2009) and log-concave distributions (Litvak et al., 2005). As such, the the above result can be easily generalized to a much larger class of random matrices.

# 6. Robustness against manipulators

In this section, we apply our results to study the "profile injection" attacks on collaborative filtering. According to the empirical study of Mobasher et al. (2006), matrix factorization, as a model-based CF algorithm, is more robust to such attacks compared to similarity-based CF algorithms such as kNN. However, as Cheng & Hurley (2010) pointed out, it may not be a conclusive argument that model-based recommendation system is robust. Rather, it may due to the fact that that common attack schemes, effective to similarity based-approach, do not exploit the vulnerability of the model-based approach.

Our discovery is in tune with both Mobasher et al. (2006) and Cheng & Hurley (2010). Specifically, we show that factorization methods are resilient to a class of common attack models, but are not so in general.

---

[4] Hence $N$ is also the orthonormal basis of any $Y$ generated with $G$ being its left multiplier.

## 6.1. Attack models

Depending on purpose, attackers may choose to inject "dummy profiles" in many ways. Models of different attack strategies are surveyed in Mobasher et al. (2007). For convenience, we propose to classify the models of attack into two distinctive categories: *Targeted Attack* and *Mass Attack*.

**Targeted Attacks** include average attack (Lam & Riedl, 2004), segment attack and bandwagon attack (Mobasher et al., 2007). The common characteristic of targeted attacks is that they *pretend* to be the honest users in all ratings except on a few targets of interest. Thus, each dummy user can be decomposed into:

$$e = e^{gnd} + s,$$

where $e^{gnd} \in \mathcal{N}$ and $s$ is sparse.

**Mass Attacks** include random attack, love-hate attack (Mobasher et al., 2007) and others. The common characteristic of mass attacks is that they insert dummy users such that many entries are manipulated. Hence, if we decompose a dummy user,

$$e = e^{gnd} + e^{gnd\perp},$$

where $e^{gnd} = \mathbb{P}^{\mathcal{N}} e$ and $e^{gnd\perp} = (I - \mathbb{P}^{\mathcal{N}})e \in \mathcal{N}^{\perp}$, then both components can have large magnitude. This is a more general model of attack.

## 6.2. Robustness analysis

By definition, injected user profiles are column-wise: each dummy user corresponds to a corrupted column in the data matrix. For notational convenience, we re-arrange the order of columns into $[\, Y \,|\, E \,]$, where $Y \in \mathbb{R}^{m \times n}$ is of all honest users, and $E \in \mathbb{R}^{m \times n_e}$ contains all dummy users. As we only care about the prediction of honest users' ratings, we can, without loss of generality, set ground truth to be $[\, Y \,|\, E^{gnd} \,]$ and the additive error to be $[\, 0 \,|\, E^{gnd\perp} \,]$. Thus, the recovery error is $Z = [\, Y^* - Y \,|\, E^* - E^{gnd} \,]$.

**Proposition 3.** *Assume all conditions of Theorem 1 hold. Under "Targeted Attacks", there exists an absolute constant $C$, such that*

$$\text{RMSE} \leq 4k\sqrt{\frac{s_{max}n_e}{|\Omega|}} + Ck\left( \frac{(n + n_e)r\log(n + n_e)}{|\Omega|} \right)^{\frac{1}{4}}. \tag{12}$$

*Here, $s_{max}$ is maximal number of targeted items of each dummy user.*

*Proof.* In the case of "Targeted Attacks", we have (recall that $k = \max_{(i,j)} |Y_{i,j}|$)

$$\|E^{gnd\perp}\|_F < \sum_{i=1,\dots,n_e} \|s_i\| \leq \sqrt{n_e s_{max}(2k)^2}.$$

Substituting this into Theorem 1 establishes the proposition. □

**Remark 1.** Proposition 3 essentially shows that matrix factorization approach is robust to the targeted attack model due to the fact that $s_{max}$ is small. Indeed, if the sampling rate $|\Omega|/(m(n+n_e))$ is fixed, then RMSE converges to zero as $m$ increases. This coincides with empirical results on Netflix data (Bell & Koren, 2007). In contrast, similarity-based algorithms (kNN) are extremely vulnerable to such attacks, due to the high similarity between dummy users and (some) honest users.

It is easy to see that the factorization method is less robust to mass attacks, simply because $\|E^{gnd\perp}\|_F$ is not sparse, and hence $s_{max}$ can be as large as $m$. Thus, the right hand side of (12) may not diminish. Nevertheless, as we show below, if the number of "Mass Attackers" does not exceed certain threshold, then the error will mainly concentrates on the $E$ block. Hence, the prediction of the honest users is still acceptable.

**Proposition 4.** *Assume sufficiently random subspace N (i.e., Propostion 2 holds), above definition of "Mass Attacks", and condition number $\kappa$. If $n_e < \frac{\sqrt{n}}{\kappa^2 r}(\frac{\mathbf{E}|Y_{i,j}|^2}{k^2})$ and $|\Omega| = pm(n+n_e)$ satisfying $p > 1/m^{1/4}$, furthermore individual sample rate of each users is bounded within $[p/2, 3p/2]$,[5] then with probability of at least $1 - cm^{-10}$, the RMSE for honest users and for manipulators satisfies:*

$$\text{RMSE}_Y \leq C_1 \kappa k \left(\frac{r^3 \log(n)}{p^3 n}\right)^{1/4}, \quad \text{RMSE}_E \leq \frac{C_2 k}{\sqrt{p}},$$

*for some universal constant $c$, $C_1$ and $C_2$.*

The proof of Proposition 4, deferred in the supplementary material, involves bounding the prediction error of each individual users with Theorem 4 and sum over $Y$ block and $E$ block separately. Subspace difference $\rho$ is bounded with Theorem 1 and Theorem 3 together. Finally, $\sigma_{min}$ is bounded via Proposition 2.

### 6.3. Simulation

To verify our robustness paradigm, we conducted simulation for both models of attacks. $Y$ is generated by multiplying two $1000 \times 10$ gaussian random matrix and $n_e$ attackers are appended to the back of $Y$. Targeted Attacks are produced by randomly choosing from a column of $Y$ and assign 2 "push" and 2 "nuke" targets to 1 and -1 respectively. Mass Attacks are generated using uniform distribution. Factorization is performed using ALS. The results of the simulation are summarized in Figure 1 and 2. Figure 1 compares the RMSE under two attack models. It shows that when the number of attackers increases, RMSE under targeted attack remains small, while RMSE under random attack significantly increases. Figure 2 compares $\text{RMSE}_E$ and $\text{RMSE}_Y$ under random attack. It shows that when $n_e$ is small, $\text{RMSE}_Y \ll \text{RMSE}_E$. However, as $n_e$ increases, $\text{RMSE}_Y$ grows and eventually is comparable to $\text{RMSE}_E$. Both figures agree with our theoretic prediction.

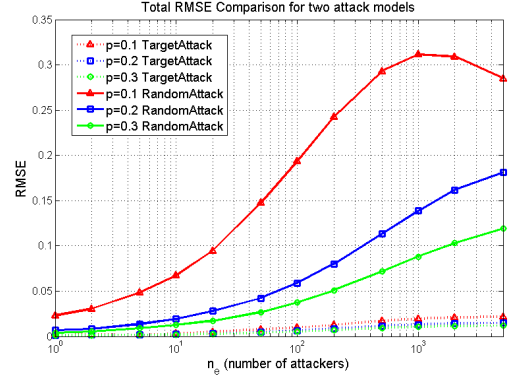

*Figure 1.* Comparison of two attack models.

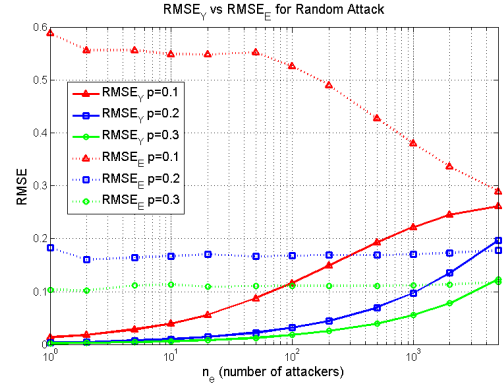

*Figure 2.* Comparison of $\text{RMSE}_Y$ and $\text{RMSE}_E$ under random attack.

## 7. Concluding discussions

This paper presented a comprehensive study of the stability of matrix factorization methods. The key results include a near-optimal stability bound, a subspace stability bound and a worst-case bound for individual columns. Then the theory is applied to the notorious manipulator problem in collaborative filtering, which leads to an interesting insight of MF's inherent robustness.

Matrix factorization is an important tool both for matrix completion task and for PCA with missing data. Yet, its practical success hinges on its stability – the ability to tolerate noise and corruption. This paper is a first attempt to understand the stability of matrix factorization, which we hope will help to guide the application of matrix factorization methods.

---

[5]This assumption is made to simplify the proof. It easily holds under i.i.d sampling.

We list some possible directions to extend this research in future. In the theoretical front, the arguably most important open question is that under what conditions matrix factorization can reach a solution near global optimal. In the algorithmic front, we showed here that matrix factorization methods can be vulnerable to general manipulators. Therefore, it is interesting to develop a robust variation of MF that provably handles arbitrary manipulators.

## Acknowledgments

This research is partially supported by the National University of Singapore under startup grant R-265-000-384-133.